\documentclass[12pt]{article}
\usepackage[affil-it]{authblk}
\usepackage[english]{babel}
\usepackage{blindtext}

\usepackage[utf8]{inputenc}
\usepackage{amsmath,amssymb,amsthm,amsfonts}
\usepackage{bm}
\usepackage[hidelinks,colorlinks = true]{hyperref}
\usepackage{enumitem}
\usepackage{blkarray}
\usepackage{tikz-cd}

\setenumerate{labelindent=0.25\parindent,leftmargin=*}
\usepackage[a4paper, total={16.5cm, 24cm}]{geometry}
\DeclareMathOperator{\rk}{rk}

\DeclareMathOperator{\Img}{Im}

\DeclareMathOperator{\codim}{codim}

\def\F{{\mathbb F}}

\def\M{{\mathcal M}}
\def\N{{\mathbb N}}
\def\SCS{{\mathbf S}}

\let\det\undefined
\DeclareMathOperator{\det}{det}

\DeclareMathOperator{\sgn}{sgn}

\newenvironment{psmallmatrix}
  {\left(\begin{smallmatrix}}
  {\end{smallmatrix}\right)}

\theoremstyle{plain}
\newtheorem{theorem}{Theorem}[section]
\newtheorem*{theorem*}{Theorem}
\newtheorem{lemma}[theorem]{Lemma}
\newtheorem*{lemma*}{Lemma}
\newtheorem{claim}[theorem]{Claim}
\newtheorem{proposition}[theorem]{Proposition}
\newtheorem*{proposition*}{Proposition}
\newtheorem{corollary}[theorem]{Corollary}
\newtheorem*{corollary*}{Corollary}

\theoremstyle{definition}
\newtheorem{definition}[theorem]{Definition}

\newtheorem{example}[theorem]{Example}

\newtheorem*{counterexample*}{Counterexample}

\newtheorem{remark}[theorem]{Remark}

\numberwithin{equation}{section}

\newcommand*{\doi}[1]{DOI: \href{http://dx.doi.org/#1}{#1}}

\begin{document}

\title{Linear maps preserving the Cullis'
determinant. III}
\author[1]{Alexander Guterman}
\author[2,*]{Andrey Yurkov}
\affil[1,2]{Department of Mathematics, Bar Ilan University, Ramat-Gan, 5290002, Israel}
\affil[*]{Corresponding author: Andrey Yurkov, andrey.yurkov@biu.ac.il}

\maketitle

\begin{abstract}
This paper is the third in the series of papers devoted to the explicit description of linear maps preserving the Cullis' determinant of rectangular matrices of size $n\times k$ with entries belonging to an arbitrary ground field which is large enough.

In this part we study the remaining case where $k = 3$ and $n \ge k + 2 = 5$ and prove that every linear map preserving the Cullis' determinant admits the same description as for the case where $k \ge 4$ and $n \ge k + 2$ obtained previously. 
\end{abstract}


\section{Introduction}\label{sec:intro}

The theory of linear maps preserving matrix invariants is a broad field of research in Linear Algebra which is related with other fields of mathematics. It attracts an attention of many scientists since the beginning of 20th century when Frobenius in~\cite{GF} obtained a description of linear maps preserving the determinant of square matrices. 

\begin{theorem}[Frobenius, {\cite[\textsection 7, Theorem I]{GF}}]\label{thm:1}
Let $S\colon \M_{n}(\mathbb C) \to \M_{n}(\mathbb C)$ be a bijective linear map satisfying $\det (S(X))=\det (X)$ for all $X \in \M_{n} (\mathbb C)$, where $\mathbb C$ denotes the field of complex numbers and $\M_{n}(\mathbb C)$  denotes the set $n\times n$ matrices with entries from $\mathbb C.$ Then there exist matrices $M, N \in \M_n(\mathbb C)$ with $\det (MN) = 1$ such that
$$S(X) = MXN\;\;\mbox{for all}\; X \in \M_{n} (\mathbb C)
\;\;\text{or}\;\;
S(X) = MX^{t}N\;\;\mbox{for all}\; X \in \M_{n}(\mathbb C).$$
\end{theorem}

This area of study continues to be a subject of active research. We refer the reader to~\cite{LAMA199233} for a survey of results up to the end of the 20th century and to~\cite{hogben_handbook_2014} for a brief introduction. In particular, there are obtained generalisations of Theorem~\ref{thm:1} regarding polynomial invariants of square matrices defined similarly including permanent~\cite{Botta1967}, immanants~\cite{DUFFNER1994567} and Schur matrix functionals~\cite{DESEGUINSPAZZIS201963}.

The Cullis' determinant is an invariant of rectangular matrix generalizing the notion of the ordinary determinant of square matrix. It was introduced by Cullis in~\cite{cullis1913} and studied later by Radi\'{c} in~\cite{radic1966}--\nocite{radic2005,radic1991}\cite{radic2008}, Makarewicz, Mozgawa, Pikuta and Sza\l{}kowski in~\cite{makarewicz2014}--\nocite{makarewicz2016}\cite{makarewicz2020}, Amiri, Fathy and Bayat in~\cite{amiri2010}, and Nakagami and Yanay in~\cite{NAKAGAMI2007422}. The Cullis' determinant has many properties similar to those of the ordinary determinant (Theorem~\ref{thm:DetNKProperties}). Together with
the Frobenius' theorem above, this motivates to characterize the linear maps preserving the Cullis' determinant. This direction of research was completely unexplored until 2024, when the authors solved the linear preserver problem for the Cullis determinant of rectangular matrices of size $n\times k$ with $n = k + 1$. As it is shown in~\cite{Guterman2024}, this case is reduced to the case $n = k$, where it is possible to use the Frobenius' theorem. 

In further studies of this problem, we obtained that for $n \ge k + 2$ and $k \ge 4$, the characterization depends on parity of $n + k$ (see~\cite{Guterman2025} and~\cite{Guterman2025b}). This result is formulated in Theorem~\ref{thm:THMNPLUSKEVEN} and Theorem~\ref{thm:THMNPLUSKODD}. Let us introduce the basic notions and definitions used throughout our paper in order to formulate these theorems. 

\vspace{\topsep}

By $\F$ we denote a field without any restrictions on its characteristic. We will use the inequality $|\F| > k$ for $k \in \mathbb N$ which means that either $\F$ is infinite or $\F$ is finite and $|\F| > k$.

We denote by $\M_{n\, k}(\F)$ the set of all $n\times k$ matrices with entries from a certain field $\F.$ Similarly, by $\M_{n\, k}(\F[x_1,\ldots, x_m])$ we denote the set of all $n\times k$ matrices with entries from the ring $\F[x_1,\ldots, x_m]$ of polynomials in variables $x_1,\ldots, x_m$ with coefficients from $\F.$ 

$O_{n\, k}\in \M_{n\, k}(\F)$ denotes the matrix with all the entries equal to zero. $I_{n} \in \M_{n\, n}(\F)$ denotes the identity matrix. Let us denote by $E_{ij} \in \M_{n\, k}(\mathbb F)$ the matrix, whose entries are all equal to zero besides the entry on the intersection of the $i$-th row and the $j$-th column, which is equal to one. By $x_{i\, j}$ we denote the element of a matrix $X$ lying on the intersection of its $i$-th row and $j$-th column. If $X \in \M_{n\,k}(\F)$, then by $\rk(X)$ we denote the rank of $X$. For a set $S$ and $i, j \in S$ by $\delta_{i\,j}$ we denote the \emph{Kronecker delta} of $i$ and $j$, which is equal to $1$ if $i = j$ and equal to $0$ otherwise.

For $A \in \M_{n\,k_1}(\mathbb F)$ and $B \in \M_{n\,k_2}(\mathbb F)$ by $A|B \in \M_{n\,k_1+k_2}(\mathbb F)$ we denote the block matrix defined by $A|B = \begin{pmatrix} A & B\end{pmatrix}$.

For $A \in \M_{n\,k}(\mathbb F)$ by $A^t \in \M_{k\,n}(\mathbb F)$ we denote the transpose of the matrix $A$, i.e. $A^{t}_{i\,j} = A_{j\,i}$ for all $1 \le i \le k, 1 \le j \le n$.

We use the notation for submatrices following~\cite{Minc1984} and~\cite{Pierce1979}. That is, by $A[J_1|J_2]$ we denote the $|J_1 |\times |J_2|$ submatrix of $A$ lying on the intersection of rows with
the indices from $J_1$ and the columns with the indices from $J_2$. By $A(J_1|J_2)$ we denote the submatrix of $A$ derived from it by striking out from it the rows with indices belonging to $J_1$ and the columns with the indices belonging to $J_2$. If the set is absent, then it means the empty set, i.e. $A(J_1|)$ denotes the matrix derived from $A$ by striking out from it the rows with indices belonging to $J_1$. We may skip curly brackets, i.e. $A[1,2|3,4] = A[\{1,2\}|\{3,4\}]$. The notation with mixed brackets is also used, i.e. $A(|1]$ denotes the first column of the matrix $A$. This notation is also used for vectors as well. In this case vectors are considered as $n\times 1$ or $1 \times n$ matrices. 

We use the bold font to denote vectors and lower indices to denote their coordinates. In the case if we need a series of vectors, we use upper indices placed in braces. For example, if $\mathbf v = \begin{pmatrix}1 \\ 0\end{pmatrix},$ then $\mathbf v^t = \begin{pmatrix}1 & 0\end{pmatrix}$ and $\mathbf v_1 = 1$. If $\mathbf u^{(1)} = \begin{pmatrix}1 \\ 0\end{pmatrix}, \mathbf u^{(2)} = \begin{pmatrix}0 \\ 1\end{pmatrix}$, then $\mathbf u^{(1)}_1 = 1$ and $\mathbf u^{(2)}_1 = 0$.\pagebreak

For $A = (a_{i\,j}) \in \M_{n\,n}(\F)$ we denote by $\det(A)$ the ordinary determinant of $A$. In addition, the following notion is used 
$$\det(A) = \begin{vmatrix}a_{1\,1} & \cdots & a_{n\,1}\\ \vdots & \ddots & \vdots\\a_{n\,1} & \cdots & a_{n\,n}\end{vmatrix}.$$
\vspace{\topsep}

\begin{definition}[{Cf.~\cite[\S3, the definition on p. 12 and \S30, Theorem~I.]{cullis1913}}]\label{def:CullisDet}
Let  $n \ge k\ge 1$ be integers, $X \in \M_{n\,k} (\mathbb F)$ be a matrix. Then Cullis' determinant $\det_{n\, k}(X)$ of $X$ is defined by
\begin{equation*}
\det_{n\,k}(X) = (-1)^{(1 + \ldots + k)}\sum_{1 \le c_1 < \ldots < c_k \le n } (-1)^{c_1 + \ldots + c_k} \begin{vmatrix}x_{c_1\, 1} & \cdots & x_{c_1\, k}\\ \vdots & \ddots & \vdots\\ x_{c_k\, 1} & \cdots & x_{c_k\, k}\end{vmatrix}.
\end{equation*}
That is, the Cullis' determinant of $X$ is an alternating sum of basic minors of $X$. We also denote $\det_{n\, k} (X)$ as follows
\[
\det_{n\,k}(X)=\begin{vmatrix}X_{1\,1} & \cdots & X_{1\,k}\\
\vdots & \cdots & \vdots\\
X_{n\,1} & \cdots & X_{n\,k}
\end{vmatrix}_{n\,k}.
\]
In the case if $n = k$, then $\det_{n\,n}$ is also denoted as $\det_{k}$ and is clearly equal to the ordinary determinant of a square matrix.

For $k > n$ the Cullis' determinant of $X \in \M_{n\,k}(\F)$ is defined by $\det_{n\, k}(X) = \det_{k\,n}(X^t)$ (see~\cite[Definition~1]{radic1966}). 
\end{definition}

\vspace{\topsep}

\begin{theorem}[{\cite[Theorem~5.14]{Guterman2025}}]\label{thm:THMNPLUSKEVEN}
Assume that $|\F| > k \ge 4, n \ge k + 2$ and  $n + k$ is even. Let $T\colon \M_{n\, k} (\F) \to \M_{n\, k} (\F)$ be a linear map. Then $\det_{n\, k} (T(X)) = \det_{n\,k}(X)$ for all $X \in \M_{n\, k} (\F)$ if and only if there exist $A \in \M_{n\, n}(\F)$ and $B \in \M_{k\, k}(\F)$ such that
\begin{equation*}
\det_{n\, k} \Bigl(A(|i_1,\ldots, i_k]\Bigr) \cdot \det_k \Bigl(B\Bigr) = (-1)^{i_1 + \ldots + i_k - 1 - \ldots - k}
\end{equation*}
for all increasing sequences $1 \le i_1 < \ldots < i_k \le n $ and
\[
T(X) = AXB\;\;\mbox{for all}\;\;X \in \M_{n\, k} (\F).
\]

\end{theorem}

\begin{theorem}[{Cf.~\cite[Theorem~4.12]{Guterman2025b}}]\label{thm:THMNPLUSKODD}
Assume that $|\F| > k \ge 4$, $n \ge k + 2$ and $n + k$ is odd. Let $T\colon \M_{n\,k} (\F) \to \M_{n\,k} (\F)$ be a linear map. Then $\det_{n\, k} (T(X)) = \det_{n\,k}(X)$ for all $X \in \M_{n\,k} (\F)$ if and only if  there exist $A \in \M_{n\, n}(\F)$ and $B \in \M_{k\, k}(\F)$ such that
\begin{equation*}
\det_{n\,k} \Bigl(A(|i_1,\ldots, i_k]\Bigr) \cdot \det_k \Bigl(B\Bigr) = (-1)^{i_1 + \ldots + i_k - 1 - \ldots - k}
\end{equation*}
for all increasing sequences $1 \le i_1 < \ldots < i_k \le n $ and a linear map $\phi\colon \M_{n\,k}(\F) \to W_{n\,k}$ such that
\begin{equation*}
T(X) = AXB + \phi(X)\;\;\mbox{for all}\;\;X \in \M_{n\,k} (\F),
\end{equation*}
 where  $W_{n\,k} \subseteq \M_{n\,k}(\F)$ denotes the space of matrices, all rows of which are equal. 
\end{theorem}

As we mentioned above, these two theorems provide a description of linear maps preserving $det_{n\,k}$ for $n \ge k + 2$ and $k \ge 4$ in the case if the ground field is large enough.

In~\cite{Guterman2025} we also consider the linear maps preserving the Cullis' determinant of $n\times 1$ and $n \times 2$ matrices, where the Cullis' determinant appears to be a linear and a quadratic form, correspondingly. The remaining case of $n \times 3$ matrices could not be carried out by the approach used in the previous papers because our argumentation there is based on the following lemma:

\begin{lemma}[{\cite[Lemma~5.8]{Guterman2025}}]Suppose that $k \ge 4$, $n \ge k$, $B \in \M_{n\, k}(\F)$ and $$\deg_{\lambda} (\det_{n\, k} (A + \lambda B)) \le 1$$ for all $A \in \M_{n\, k}(\F).$ Then $\rk (B) \le 1.$
\end{lemma}  

We use it  to show that every linear map preserving $\det_{n\,k}$ should preserve the set of matrices of rank $1$ and this lemma is not true for $k = 3$ as it is shown by the example below:

\begin{example}[{\cite[Example~5.11]{Guterman2025}}]\label{cntex:RK2DetN3LowDeg}Suppose that $B = \begin{psmallmatrix}
                                 1 & 0 & 0 & 0 & -1 & 0 & \cdots & 0\\
                                 0 & 1 & 0 & -1 & 0 & 0 & \cdots & 0\\
                                 0 & 0 & 0 & 0 & 0 & 0 & \cdots & 0
                                        \end{psmallmatrix}^t \in \M_{n\,3}(\F).$ Then $\rk (B) = 2$ and $\deg_{\lambda} (\det_{n\,k}(A + \lambda B)) \le 1$ for all $A \in \M_{n\,3}(\F)$.
\end{example}

To overcome this obstacle in the case if $k = 3$ and $n \ge 5 = 3 + 2$, we prove in Lemma~\ref{lem:N3LowDegImpliesBInU} that every matrix $B \in \M_{n\,3}(\F)$ with $\rk (B) = 2$ such that $\deg_{\lambda} \left(\det_{n\,3}(A + \lambda B)\right) \le 1$ belongs to the space of matrices whose alternating sum of rows is the zero vector. Then, using certain kind of ``density argument'' (Lemma~\ref{lem:Rk1NonzeroRk1All}), we prove that  every linear map preserving $\det_{n\,3}$ preserves the set of matrices of rank $1$ which allows us to conclude that the set of linear maps preserving $\det_{n\,3}$ has the same description as for the case $k > 3$ discussed above and prove Theorem~\ref{thm:MainTheoremEvenINCLK3} and Theorem~\ref{thm:MainTheoremOddINCLK3} below. Accordingly, the main result of this paper could be formulated as follows.

\begin{proposition}[Theorem~\ref{thm:MainTheoremEvenINCLK3}, Theorem~\ref{thm:MainTheoremOddINCLK3}]The statements of Theorem~\ref{thm:THMNPLUSKEVEN} and Theorem~\ref{thm:THMNPLUSKODD} hold also for $k = 3$ and $n \ge k + 2 = 5$. 
\end{proposition}

Thus, the result of Frobenius, together with our results, provides a complete solution to the problem of describing linear maps preserving the Cullis determinant in the case where the ground field is sufficiently large. We conjecture that the assumption on the size of the ground field could be omitted, although a different approach would then be required.

The further text is organized as follows: in the sections~\ref{sec:prelimalglinalg} and~\ref{sec:prelimcul} we provide the preliminary material from Algebra, Linear Algebra and the theory of the Cullis' determinant; in Section~\ref{sec:LinPresDetN3} we obtain the description of linear maps preserving $\det_{n\,3}$ for $n \ge 5$.

\section{Preliminaries from Algebra and Linear Algebra}\label{sec:prelimalglinalg}

First recall the following general facts regarding polynomials over a field which are used in the proof of Lemma~\ref{lem:Rk1NonzeroRk1All}.

\begin{lemma}[{Cf.~\cite[IV, \S 1, Corollary 1.7]{Lang2002}}]Let $\F$ be an infinite field. Let $f$ be a
polynomial in $n$ variables over $\F$. If $f$ defines the zero function on $\F^{n}$, then $f = 0$.
\end{lemma}

\begin{lemma}[{Cf.~\cite[IV, \S 1, Corollary 1.8]{Lang2002}}]Let $\F$ be a finite field with $q$ elements. Let $f$ be a
polynomial in $n$ variables over $\F$ such that the degree of $f$ in each variable
is less than $q$. If $f$ defines the zero function on $\F^{n}$, then $f = 0$.
\end{lemma}

\begin{corollary}\label{DegFLessZero}Let $\F$ be an arbitrary field and $f$ be a
polynomial in $n$ variables over $\F$ such that either $\F$ is infinite or the degree of $f$ in each variable
is less than $|\F|$. If $f(\mathbf x) = 0$ for all $\mathbf x \in \F^n$, then $f = 0$.
\end{corollary}

The ``density'' argument proved in the following lemma was mentioned in the introduction. It is used to 
deduce  Lemma~\ref{lem:LinPresRankK3} from Lemma~\ref{cor:RkAEq1NotInU3ImpliesRkTAEq1}.

\begin{lemma}\label{lem:Rk1NonzeroRk1All}Assume that $|\F| > 3$ and $n \ge k \ge 1$, $m \ge 1$ are integers. Let $c_1, \ldots, c_m \in \F$ and $B^{(1)}, \ldots, B^{(m)} \in \M_{n\,k}(\F)$ be such that $(c_1, \ldots, c_m) \neq \mathbf 0$. If 
\begin{equation}\label{lem:Rk1NonzeroRk1All:eq}
\rk\left(\alpha_1 B^{(1)} + \ldots + \alpha_n B^{(m)}\right) \le 1\;\;\;\;\mbox{for all}\;\;\alpha_1, \ldots, \alpha_m \in \F\;\;\mbox{such that}\;\;\sum\limits_{i=1}^{m} c_i\alpha_i\neq 0,
\end{equation}
then
\begin{equation}\label{lem:Rk1NonzeroRk1All:eqq}
\rk\left(\alpha_1 B^{(1)} + \ldots + \alpha_n B^{(m)}\right) \le 1\;\;\;\;\mbox{for all}\;\;\alpha_1, \ldots, \alpha_m \in \F\;\;\mbox{such that}\;\; (\alpha_1, \ldots, \alpha_m) \neq \mathbf 0.
\end{equation}
\end{lemma}
\begin{proof}
By contradiction. Suppose that  $B^{(1)}, \ldots, B^{(m)} \in \M_{n\,k}(\F)$ satisfy the condition~\eqref{lem:Rk1NonzeroRk1All:eq} and there exists $\alpha_1, \ldots, \alpha_n \in \F$ such that $(\alpha_1, \ldots, \alpha_m) \neq \mathbf 0$ and 
\begin{equation}\label{lem:Rk1NonzeroRk1All:eq101}
\rk\left(A\right) > 1,
\end{equation}
where $A \in \M_{n\,k}(\F)$ is defined by 
\begin{equation}\label{lem:Rk1NonzeroRk1All:eq102}
A = \alpha_1 B^{(1)} + \ldots + \alpha_m B^{(m)}.
\end{equation}
Therefore, there exist  $1 \le i_1, i_2 \le n$ and $1 \le j_1, j_2 \le k$ satisfying
\begin{equation}\label{lem:Rk1NonzeroRk1All:eq1}
\begin{vmatrix}a_{i_1\, j_1} & a_{i_1\,j_2}\\ a_{i_2\,j_1} & a_{i_2\,j_2}\end{vmatrix} \neq 0.
\end{equation}
Note here that $k \ge 2$ because if $k = 1$, then every matrix in $\M_{n\,k}(\F)$ has rank at most one which contradicts with~\eqref{lem:Rk1NonzeroRk1All:eq101}.

Let $Q(t_1,\ldots, t_m) \in \F[t_1,\ldots, t_m]$ be a polynomial defined by
\begin{equation}\label{lem:Rk1NonzeroRk1All:eq10005}
Q(t_1,\ldots, t_m) = \begin{vmatrix}a_{i_1\, j_1}(t_1,\ldots, t_m) & a_{i_1\,j_2}(t_1,\ldots, t_m)\\ a_{i_2\,j_1}(t_1,\ldots, t_m) & a_{i_2\,j_2}(t_1,\ldots, t_m)\end{vmatrix},
\end{equation}
where $A(t_1,\ldots, t_m) = \left(a_{i\,j}(t_1,\ldots, t_m)\right) \in \M_{n\,k}\left(\F[t_1,\ldots, t_m]\right)$ is defined by
$$A(t_1,\ldots, t_m) = t_1 B^{(1)} + \ldots + t_m B^{(m)}.$$
Since 
$Q(\alpha_1, \ldots, \alpha_m) = \det\big(A[i_1,i_2|j_1,j_2]\big) \neq 0$
by~\eqref{lem:Rk1NonzeroRk1All:eq1}, then 
\begin{equation}\label{lem:Rk1NonzeroRk1All:eq10001}
Q\neq 0
\end{equation}
as an element of $\F[t_1,\ldots, t_m]$. In addition, 
\begin{equation}\label{lem:Rk1NonzeroRk1All:eq10002}
\deg(Q) \leq 2
\end{equation} by the definition of $Q$. 

Consider a polynomial $P(t_1,\ldots, t_m) \in \F[t_1,\ldots, t_m]$ defined by
\begin{equation}\label{lem:Rk1NonzeroRk1All:eq105}
 P = \left(\sum\limits_{i=1}^{m}c_i t_i\right)\cdot Q.
\end{equation}
The equality~\eqref{lem:Rk1NonzeroRk1All:eq10001} implies that $P$ is nonzero because $(c_1, \ldots, c_m) \neq \mathbf 0$ by the assumption of the lemma. In addition, 
\begin{equation}\label{lem:Rk1NonzeroRk1All:eq10007}
\deg(P) = \deg\left(\sum\limits_{i=1}^{m}c_i t_i\right) + \deg(Q) \le 1 + 2 = 3
\end{equation}
by~\eqref{lem:Rk1NonzeroRk1All:eq10002}. Since $|\F| > 3$ by the conditions of the lemma, then $\deg(P) < |\F|$ by~\eqref{lem:Rk1NonzeroRk1All:eq10007}. Therefore, by Corollary~\ref{DegFLessZero} there exist $\alpha'_1, \ldots, \alpha'_m \in \F$ such that 
\begin{equation}\label{lem:Rk1NonzeroRk1All:eq104}
P(\alpha'_1, \ldots, \alpha'_m) \neq 0.
\end{equation}
Let $A' \in \M_{n\,k}(\F)$ be defined by
$$A' = \alpha'_1 B^{(1)} + \ldots + \alpha'_m B^{(m)}.$$
The definition of $P$ in~\eqref{lem:Rk1NonzeroRk1All:eq105} together with~\eqref{lem:Rk1NonzeroRk1All:eq104} imply
 that 
 \begin{equation}\label{lem:Rk1NonzeroRk1All:eq10003}
 \sum_{i=1}^{m} c_i \alpha'_i \neq 0
 \end{equation}
 and
 \begin{equation}\label{lem:Rk1NonzeroRk1All:eq10004}
 Q(\alpha'_1, \ldots, \alpha'_m) \neq 0.
 \end{equation}
 Recall that $Q(\alpha'_1, \ldots, \alpha'_m)$ is equal to $\det\big(A'[i_1,i_2|j_1,j_2]\big)$ by the definition of $Q$ in~\eqref{lem:Rk1NonzeroRk1All:eq10005}. Hence, the equality~\eqref{lem:Rk1NonzeroRk1All:eq10004} implies that
 \begin{equation}\label{lem:Rk1NonzeroRk1All:eq10006}
 \rk(A') > 1.
 \end{equation}

Given the equalities~\eqref{lem:Rk1NonzeroRk1All:eq10003} and~\eqref{lem:Rk1NonzeroRk1All:eq10006}, we reach a contradiction with~\eqref{lem:Rk1NonzeroRk1All:eq}. Thus, $B^{(1)}, \ldots, B^{(m)}$ satisfy the condition~\eqref{lem:Rk1NonzeroRk1All:eqq}.
 \end{proof}

The following theorem due to Westwick is reformulated in order to apply it to maps on matrix spaces. This theorem is used in the final step of the proof of Theorem~\ref{thm:MainTheoremEvenINCLK3}.

\begin{lemma}[\cite{westwick1967}, Cf. Theorem 3.5]\label{lem:westwick}Assume that $n > k$ and let $T$ be an onto linear map of $\M_{n\,k}(\F)$ such that $\rk(T(X)) = 1$ for all $X \in \M_{n\,k}(\F)$ with $\rk(X) = 1$. Then there exist invertible matrices $A \in \M_{n\,n}(\F)$ and $B \in \M_{k\,k}(\F)$ such that
  \[
    T(X) = AXB\;\;\mbox{for all}\;\; X \in \M_{n\,k}(\F).
  \]
\end{lemma}\pagebreak

\section{Preliminaries from the theory of the Cullis' determinant}\label{sec:prelimcul}

In this section we present the definitions and propositions from the theory of the Cullis' determinant used in the sequel.

We begin with the list of the properties of $\det_{n\, k}$ which are similar to corresponding properties of the ordinary determinant (see \cite[\textsection 5, \textsection 27, \textsection 32]{cullis1913}, \cite{NAKAGAMI2007422} and~\cite{radic1966} for detailed proofs).

\begin{theorem}\label{thm:DetNKProperties}
\begin{enumerate}
\item[]
\item For $X \in \M_{n}(\mathbb F),$ $\det_{n\,n}(X) = \det (X).$
\item For $X \in \M_{n\,k}(\mathbb F),$ $\det_{n\,k}(X)$ is a linear function of columns of $X$.
\item If a matrix $X \in \M_{n\,k}(\mathbb F)$ has two identical columns or one of its columns is a linear combination of other columns, then $\det_{n\,k}(X)$ is equal to zero.
\item For $X \in \M_{n\,k}(\mathbb F),$ interchanging any two columns of $X$ changes the sign of $\det_{n\,k}(X)$.
\item Adding a linear combination of columns of $X$ to another column of $X$ does not change $\det_{n\,k}(X)$.
\item For $X \in \M_{n\,k}(\mathbb F),$ $\det_{n\,k}(X)$ can be calculated using the Laplace expansion along a column of $X$ (see Lemma~\ref{lem:DetNKLaplaceExp} for the precise formulation).
\end{enumerate}
\end{theorem}

\begin{definition}\label{def:RAB}Let $A, B \in \M_{n\,k}(\F)$ and $1 \le i_1 < \ldots < i_d \le k$. By $\mathcal R(A|\{i_1, \ldots, i_d\}|B) \in \M_{n\,k}(\F)$ we denote a matrix obtained from $A$ by replacing each column $A(|i_m]$ with the column $B(|i_m]$ for all $1 \le m \le d$.
\end{definition}

\begin{corollary}\label{cor:CullisBinomialExpansion}Let $n \ge k$, $A, B \in \M_{n\,k} (\F).$
Then
\begin{equation}\label{eq:CullisBinomialExpansion}
\det_{n\,k} (A + \lambda B)\\
= \sum_{d = 0}^{k}\lambda^d \left( \sum_{1 \le i_1 < \ldots < i_d \le k} \det_{n\,k}\Bigl(\mathcal R(A|\{i_1, \ldots, i_d\}|B)\Bigr)\right)\;\;\mbox{in}\;\; \F[\lambda].
\end{equation}
\end{corollary}
\begin{proof}This is a direct consequence from the multilinearity of $\det_{n\,k}$ with respect to the columns of a matrix.
\end{proof}
\begin{corollary}\label{cor:DegDetABLEQK}If $A, B \in \M_{n\,k}(\F)$, then $\deg_\lambda\left(\det_{n\,k} (A + \lambda B)\right) \le k$.
\end{corollary}

\begin{lemma}[{Cf.~\cite[Theorem~16]{NAKAGAMI2007422}}]\label{lem:DetNKLaplaceExp}
Let $1 < k \le n$. For any $n \times k$ matrix $X = (x_{i\, j})$ the expansion of $\det_{n\, k}(X)$ along the $j$-th column
is given by
\[
\det_{n\, k}(X) = \sum_{i = 1}^n (-1)^{i+j} x_{i\, j} \det_{(n-1)\,(k-1)}\Bigl(X(i|j)\Bigr),
\]
where $X(i|j)$ denotes the $(n-1) \times (k-1)$ matrix obtained from $Y$ by deleting the $i$-th row and the $j$-th column.
\end{lemma}

\begin{lemma}[{Cf.~\cite[Lemma 3.2]{amiri2010}}]\label{lem:AMIRILEMMA32} Let $n \ge k \ge 1$ and $x_{i\,j} \in \F$ for all $1 \le i \le n, 1 \le j \le k$. Then
\[
\begin{vmatrix}
x_{1\,1} & \ldots & x_{1\,k}\\
\vdots & \vdots & \vdots\\
x_{n\,1} & \ldots & x_{n\,k}\\
0 & \ldots & 0
\end{vmatrix}_{n+1\,k} =
\begin{vmatrix}
x_{1\,1} & \ldots & x_{1\,k}\\
\vdots & \vdots & \vdots\\
x_{n\,1} & \ldots & x_{n\,k}
\end{vmatrix}_{n\,k}.
\]
\end{lemma}

\begin{lemma}[Cf.~{\cite[Lemma~20]{NAKAGAMI2007422}}]~\label{lem:NAKAGAMILEMMA20} Assume that $n > k \ge 1$. Let  $X \in \M_{n\,k}(\F)$ and $Y \in \M_{n\,(k+1)}(\F)$ is defined by $Y = X | \begin{psmallmatrix}1 \\ \vdots \\ 1\end{psmallmatrix}$. Then
\[
\det_{n\,(k+1)}(Y) = \begin{cases}
\det_{n\,k}(X), & \mbox{$n + k$ is odd},\\
0, & \mbox{$n + k$ is even}.
\end{cases}
\]
\end{lemma}

\begin{lemma}[Cf.~{\cite[Lemma~3.1]{Guterman2025}}]\label{lem:TwoSidedMulPreservesDet}Let $n \ge k$, $A \in \M_{n\,n}(\F)$, $B \in \M_{k\,k}(\F)$, and let $T\colon \M_{n\,k}(\F) \to \M_{n\,k}(\F)$ be a linear map defined by
\begin{equation*}
T(X) = AXB\;\;\mbox{for all}\;\; X \in \M_{n\,k}(\F).
\end{equation*}
 Then
\[
\det_{n\,k} (T(X)) = \det_{n\,k}(X)\;\;\mbox{for all}\;\; X \in \M_{n\,k} (\F).
\]
 if and only if
\begin{equation}\label{eq:TwoSideMul}
\det_{n\,k} \Bigl(A(|i_1,\ldots, i_k]\Bigr)\cdot \det_k \Bigl(B\Bigr) = (-1)^{i_1 + \ldots + i_k - 1 - \ldots - k}
\end{equation}
for all increasing sequences $1 \le i_1 < \ldots < i_k \le n $. 
\end{lemma}

\begin{lemma}[Cf.~{\cite[Proposition~24]{NAKAGAMI2007422}}]\label{cor:rightmatrixmult}If $X \in \M_{n\,k}(\F)$ and $Y \in \M_{k\,k}(\F)$, then
\[
\det_{n\,k}(XY) =  \det_{n\,k}(X)\det_k(Y).
\]
\end{lemma}

\begin{lemma}[{\cite[Corollary~4.7]{Guterman2025}}]\label{lem:LinPresIso}Assume that $n \ge k + 2,$ $|\F| > k$ and $n + k$ is even. Let $T$ be a linear map on $\M_{n\,k}(\F)$ such that $\det_{n\,k}(T(X)) = \det_{n\,k}(X)$ for all $X \in \M_{n\,k}(\F)$. Then $T$ is invertible.
\end{lemma}

In our further proofs we will use the properties of linear maps $\SCS_{a\,b}$ preserving the Cullis' determinant which are defined below.

\begin{definition}[Cf.~{\cite[Definition~4.3]{Guterman2025}}]\label{def:SCSDef}Let $n \ge k$, $1 \le a \le n$ and  $1 \le b \le k$. By $\SCS_{a\,b}$ we denote a linear map on  $\M_{n\,k}(\F)$ defined by
\begin{equation*}
\SCS_{a\,b}\begin{pmatrix}
x_{1\,1} & \cdots  & x_{1\,k}\\
  \vdots & \ddots & \vdots\\
x_{n\,1} & \cdots & x_{n\,k}\\
\end{pmatrix} =(-1)^{n-a}\cdot  
\begin{pmatrix}
(-1)^{1 - \delta_{1\,b}}x_{a\,b} & \cdots & x_{a\,1} & \cdots & x_{a\,k}\\
\vdots & \ddots & \vdots & \ddots & \vdots\\
(-1)^{1 - \delta_{1\,b}}x_{n\,b} & \cdots & x_{n\,1} & \cdots & x_{n\,k}\\
-(-1)^{1 - \delta_{1\,b}}x_{1\,b} & \cdots & -x_{1\,1} & \cdots & -x_{1\,k}\\
  \vdots & \ddots & \vdots & \ddots & \vdots\\
-(-1)^{1 - \delta_{1\,b}}x_{(a-1)\,b} & \cdots & -x_{(a-1)\,1} & \cdots &  -x_{(a-1)\,k}\\
\end{pmatrix}
\end{equation*}
for all $X = (x_{i\,j}) \in \M_{n\,k}(\F)$. That is, $\SCS_{a\,b}(X)$ is obtained from $X$ by performing the following sequence of operations:
\begin{enumerate}
\item\label{lst:SIJSeq:it1}the row cyclical shift sending $a$-th row of $X$ to the first row of the result;
\item multiplying the bottom $a-1$ rows by $-1$;
\item exchanging the first and the $b$-th column;
\item\label{lst:SIJSeq:it4} multiplying the first column  by $(-1)^{1 - \delta_{1\,b}}$;
\item multiplying all the entries by $(-1)^{n-a}$.
\end{enumerate}
\end{definition}

Lemma~\ref{lem:ReduceTo11ByShift} and Lemma~\ref{lem:SemiCyclicUKImp} show that $\SCS_{a\,b}$ has the required properties.

\begin{lemma}[Cf.~{\cite[Lemma~4.5]{Guterman2025}}]\label{lem:ReduceTo11ByShift}Assume that $n \ge k$ and $n + k$ be even. Then $\SCS_{a\,b}$ is an invertible linear map preserving $\det_{n\,k}$ for all $1 \le a \le n$, $1 \le b \le k$.
\end{lemma}

\begin{lemma}\label{lem:SemiCyclicUKImp}Assume that $n \in \N$ is odd. Let $1 \le a \le n, 1 \le b \le k$, $X = (x_{i\,j}) \in \M_{n\,k}(\F)$ and $Y = (y_{i\,j})$ be defined by $Y = \SCS_{a\,b}(X)$. If $\sum_{i=1}^{n} (-1)^{i} y_{i\,1} = 0$, then $\sum_{i=1}^{n} (-1)^{i} x_{i\,b} = 0$.
\end{lemma}
\begin{proof}
Assume that
\begin{equation}\label{lem:SemiCyclicUKImp:eq102}
\sum_{i=1}^{n} (-1)^{i} y_{i\,1} = 0.
\end{equation}
Rewrite the sum $\sum_{i=1}^{n} (-1)^{i} x_{i\,b}$ as follows
\begin{equation}\label{lem:SemiCyclicUKImp:eq1003}
\sum_{i=1}^n (-1)^{i} x_{i\,b} = \sum_{i=1}^{a-1} (-1)^{i} x_{i\,b} + \sum_{i=a}^n (-1)^{i} x_{i\,b}.
\end{equation}
Using the definition of $\SCS_{a\,b}$ we obtain that
\begin{multline}\label{lem:SemiCyclicUKImp:eq100001}
\sum_{i=n-a+2}^{n} (-1)^{i - n + a - 1} y_{i\,1} = \sum_{i=n-a+2}^{n} (-1)^{i - n + a - 1} \left(\SCS_{a\,b}(X)\right)_{i\,1}\\
  = \sum_{i=1}^{a-1} (-1)^{i - n + a - 1} \left(-(-1)^{n - a + 1 - \delta_{1\, b}} x_{i\, b}\right) =  -\left((-1)^{1 - \delta_{1\,b}}\sum_{i=1}^{a-1} (-1)^{i} x_{i\,b}\right)
\end{multline}
and
\begin{multline}\label{lem:SemiCyclicUKImp:eq100002}
\sum_{i=1}^{n-a+1} (-1)^{i + a - 1}y_{i\,1} = \sum_{i=1}^{n-a+1} (-1)^{i + a - 1} \left(\SCS_{a\,b}(X)\right)_{i\,1}\\
  = \sum_{i=a}^n (-1)^{i + a - 1} \left((-1)^{n - a + 1 - \delta_{1\, b}} x_{i\, b}\right) =  (-1)^{1 - \delta_{1\,b}} \sum_{i=a}^n (-1)^{i}x_{i\,b}.
\end{multline}
Hence,
\begin{equation}\label{lem:SemiCyclicUKImp:eq1001}
\sum_{i=1}^{a-1} (-1)^{i} x_{i\,b} \overset{\eqref{lem:SemiCyclicUKImp:eq100001}}{=\joinrel=\joinrel=\joinrel=} -(-1)^{1 - \delta_{1\,b}}\left(\sum_{i=n-a+2}^{n} (-1)^{i - n + a - 1} y_{i\,1}\right)
\end{equation}
and
\begin{equation}\label{lem:SemiCyclicUKImp:eq1002}
\sum_{i=a}^n (-1)^{i} x_{i\,b} \overset{\eqref{lem:SemiCyclicUKImp:eq100002}}{=\joinrel=\joinrel=\joinrel=} (-1)^{1 - \delta_{1\,b}}\sum_{i=1}^{n-a+1} (-1)^{i + a - 1}y_{i\,1}.
\end{equation}
In addition, since $n$ is odd, then $(-1)^{-n} = -1$ and consequently
\begin{equation}\label{lem:SemiCyclicUKImp:eq1004}
-\left(\sum_{i=n-a+2}^{n} (-1)^{i - n + a - 1} y_{i\,1}\right) =\sum_{i=n-a+2}^{n} (-1)^{i + a - 1}y_{i\,1}.
\end{equation}
Therefore,
\begin{multline}\label{lem:SemiCyclicUKImp:eq101}
\sum_{i=1}^n (-1)^{i} x_{i\,b} \overset{\eqref{lem:SemiCyclicUKImp:eq1003}}{=\joinrel=\joinrel=} \sum_{i=1}^{a-1} (-1)^{i} x_{i\,b} + \sum_{i=a}^n (-1)^{i} x_{i\,b}\\
\phantom{XX}\overset{\eqref{lem:SemiCyclicUKImp:eq1001}\;\;\scriptsize\mbox{and}\;\;\eqref{lem:SemiCyclicUKImp:eq1002}}{=\joinrel=\joinrel=\joinrel=\joinrel=\joinrel=\joinrel=\joinrel=}- (-1)^{1 - \delta_{1\,b}}\left(\sum_{i=n-a+2}^{n} (-1)^{i - n + a - 1} y_{i\,1}\right) +  (-1)^{1 - \delta_{1\,b}}\sum_{i=1}^{n-a+1} (-1)^{i + a - 1} y_{i\,1}\phantom{XXX}\\
\overset{\eqref{lem:SemiCyclicUKImp:eq1004}}{=\joinrel=\joinrel=} (-1)^{1 - \delta_{1\,b}}\sum_{i=n-a+2}^{n} (-1)^{i + a - 1}y_{i\,1} + (-1)^{1 - \delta_{1\,b}}\sum_{i=1}^{n-a+1} (-1)^{i + a - 1}y_{l\,1}\\ = (-1)^{a - \delta_{1\,b}} \sum_{i=1}^n (-1)^{i}y_{i\,1}.
\end{multline}
Thus, equality~\eqref{lem:SemiCyclicUKImp:eq101} combined with the equality~\eqref{lem:SemiCyclicUKImp:eq102} implies that $\sum_{i=1}^{n} (-1)^{i} x_{i\,b} = 0.$
\end{proof}

The proofs in the next section also rely on the properties of the vector space $U_{n\,k} \subseteq  \M_{n\,k}(\F)$ defined below.

\begin{definition}\label{def:UNK}By $U_{n\,k} \subseteq \M_{n\,k}(\F)$ we denote a vector space of all matrices whose alternating row sum is the zero vector. That is,
\[
U_{n\,k} = \{X \in \M_{n\,k}(\F)\,\mid\,\sum_{i=1}^{n} (-1)^iX[i|) = \mathbf 0\}.
\]
Note that $\dim\left(U_{n\,k}\right) = (n-1)k$ and $\codim\left(U_{n\,k}\right) = k$.
\end{definition}

\begin{lemma}[{Cf.~\cite[Theorem~7.11]{Guterman2025l}}]\label{lem:MaxDimKEvenAltSumZero}
Let $n \ge k + 2$ and $V \subseteq \M_{n\,k}(\mathbb F)$ be a vector space. Then $\det_{n\,k} (X) = 0$ for all $X \in V$ and $\codim(V) = k$ if and only if $k$ is odd and $V = U_{n\,k}.$
\end{lemma}

\begin{lemma}\label{lem:DNKPreserverPreservesUn3}Let $k \ge 1$ and $n \ge k + 2$ be odd integers, $\F$ be a field with $|\F|> k$, $T\colon \M_{n\,k} (\F) \to \M_{n\,k} (\F)$ be a linear map. If $\det_{n\, k} (T(X)) = \det_{n\,k}(X)$ for all $X \in \M_{n\,k} (\F),$ then $T(U_{n\, k}) = U_{n\,k}$.
\end{lemma}
\begin{proof}
It follows from Lemma~\ref{lem:LinPresIso} that $T$ is invertible. Hence, $T$ preserves dimensions and codimensions of vector spaces. In particular,
\begin{equation}\label{lem:DNKPreserverPreservesUn3:eq1}
\codim(T(U_{n\, k})) = \codim(U_{n\, k}) = k.
\end{equation}

Let $Y \in T(U_{n\,k})$. Since $T$ preserves $\det_{n\,k}$, then 
\begin{equation}\label{lem:DNKPreserverPreservesUn3:eq2}
\det_{n\,k}(Y) = \det_{n\,k}(T(T^{-1}(Y))) = \det_{n\,k}(T^{-1}(Y)).
\end{equation}
In addition,
\begin{equation}\label{lem:DNKPreserverPreservesUn3:eq3}
\det_{n\,k}(T^{-1}(Y)) = 0
\end{equation}
by Lemma~\ref{lem:MaxDimKEvenAltSumZero} because $T^{-1}(Y) \in U_{n\,k}$. By combining~\eqref{lem:DNKPreserverPreservesUn3:eq2} and~\eqref{lem:DNKPreserverPreservesUn3:eq3}, we obtain that
\begin{equation*}
\det_{n\,k}(Y) = 0.
\end{equation*}

Thus, $\det_{n\,k}(Y) = 0$ for all $Y \in T(U_{n\,k})$ and $\codim(T(U_{n\,k})) = k$ by~\eqref{lem:DNKPreserverPreservesUn3:eq1}. Hence, the vector space $V = T(U_{n\,k})$ satisfies the conditions of Lemma~\ref{lem:MaxDimKEvenAltSumZero} and consequently $T(U_{n\, k}) = U_{n\,k}$.
\end{proof}

Similarly to~\cite{Guterman2025}, we use the correspondence between the property ``$\rk(B) \le 1$'' and the property ``$\deg_{\lambda} (\det_{n\, k} (A + \lambda B)) \le 1\;\;\mbox{for all}\;\;A \in \M_{n\, k}(\F)$'' as it is stated in the next two lemmas.

\begin{lemma}[{\cite[Lemma~5.1]{Guterman2025}}]\label{lem:Rank1CullisDeg1}Let $n \ge k$  and $B \in \M_{n\,k}(\F)$ be such that $\rk (B) \le 1.$ Then 
$$\deg_{\lambda} (\det_{n\, k} (A + \lambda B)) \le 1\;\;\mbox{for all}\;\;A \in \M_{n\, k}(\F).$$
\end{lemma}

\begin{lemma}[{\cite[Lemma~5.9]{Guterman2025}}]\label{lem:DetNKPresPresLowDeg}Assume that $n \ge k + 2,$ $|\F| > k$ and $n + k$ is even. Let $B \in \M_{n\,k}(\F)$ and $T \colon \M_{n\,k}(\F) \to \M_{n\,k}(\F)$ be a linear map such that $\det_{n\,k}(T(X)) = \det_{n\,k}(X)$ for all $X \in \M_{n\,k}(\F)$. Then
\begin{equation*}\label{lem:DetNKPresPresLowDeg:eq1}
\deg_{\lambda}(\det_{n\, k} (A + \lambda B)) \le 1\;\;\mbox{for all}\;\;A \in \M_{n\, k}(\F)
\end{equation*}
implies that
\begin{equation*}\label{lem:DetNKPresPresLowDeg:eq2}
\deg_{\lambda}(\det_{n\, k} (A + \lambda T(B))) \le 1\;\;\mbox{for all}\;\;A \in \M_{n\, k}(\F).
\end{equation*}
\end{lemma}

\section{Linear preservers of $\det_{n\,3}$}\label{sec:LinPresDetN3}
In this section we prove that linear $\det_{n\,3}$-preservers admit the same description as we obtained for $\det_{n\,k}$ when  $k \ge 4$ in~\cite{Guterman2025, Guterman2025b} (Theorem~\ref{thm:MainTheoremEvenINCLK3} and Theorem~\ref{thm:MainTheoremOddINCLK3}). For this, we establish that for odd $n$ every linear map preserving $\det_{n\,3}$ preserves the set of matrices of rank one (Lemma~\ref{lem:LinPresRankK3}). Having this fact, we prove Theorem~\ref{thm:MainTheoremEvenINCLK3} and Theorem~\ref{thm:MainTheoremOddINCLK3} using the same arguments as in~\cite{Guterman2025, Guterman2025b}. Despite this, we provide the full and explicit proofs for the sake of clarity.

In order to prove Lemma~\ref{lem:LinPresRankK3}, we show that every matrix of size $n\times 3$ which does not allow us to apply the arguments from~\cite{Guterman2025} (i.e., the matrix $B$ in Example~\ref{cntex:RK2DetN3LowDeg}) should belong to $U_{n\,3}$ (Definition~\ref{def:UNK}). Before we prove this fact, we need to establish the following equality.

\begin{lemma}\label{lem:CullisN31XY}
Let $n \ge 3$ be an odd integer. Then
\begin{equation}\label{lem:CullisN31XY:eq}
\begin{vmatrix}x_{1} & y_{1} & 1\\ \vdots & \vdots & \vdots \\ x_{n-1} & y_{n-1} & 1\\ 0 & y_{n} & 0\end{vmatrix}_{n\,3} = y_{n}\left(\sum_{i=1}^{n-1} (-1)^{i}x_i\right)\;\;\mbox{for all}\;\;x_1, \ldots, x_{n-1}, y_1, \ldots, y_n \in \F.
\end{equation}

\end{lemma}

\begin{proof}Using the Laplace expansion along the second column (Lemma~\ref{lem:DetNKLaplaceExp}) we obtain that
\begin{equation}\label{eq:CullisN31XY:1}
\begin{vmatrix}x_{1} & y_{1} & 1\\ \vdots & \vdots & \vdots \\ x_{n-1} & y_{n-1} & 1\\ 0 & y_{n} & 0\end{vmatrix}_{n\,3}
 = \left(\sum_{i=1}^{n-1}(-1)^{i+2}y_i \begin{vmatrix}x_{1} & 1\\ \vdots & \vdots \\ x_{i-1} & 1\\ x_{i+1} & 1 \\ \vdots & \vdots\\ x_{n-1} & 1\\ 0 & 0\end{vmatrix}_{(n-1)\,2}\right) + (-1)^{n+2}y_{n}\begin{vmatrix}x_{1} & 1\\ \vdots & \vdots \\\vdots & \vdots \\\vdots & \vdots \\\vdots & \vdots \\  x_{n-1} & 1\end{vmatrix}_{(n-1)\,2}.
\end{equation}

Consider the first summand on the right-hand side of~\eqref{eq:CullisN31XY:1}. It is expressed as a sum of $n-1$ terms of the form 
\begin{equation*}
(-1)^{i+2}y_i \begin{vmatrix}x_{1} & 1\\ \vdots & \vdots \\ x_{i-1} & 1\\ x_{i+1} & 1 \\ \vdots & \vdots\\ x_{n-1} & 1\\ 0 & 0\end{vmatrix}_{(n-1)\,2},\quad 1 \le i \le n - 1.
\end{equation*}
Let us show that the third multiplier in each term above is zero. Indeed, 
\begin{equation}\label{lem:CullisN31XY:eq105}
\begin{vmatrix}x_{1} & 1\\ \vdots & \vdots \\ x_{i-1} & 1\\ x_{i+1} & 1 \\ \vdots & \vdots\\ x_{n-1} & 1\\ 0 & 0\end{vmatrix}_{(n-1)\,2} =
\begin{vmatrix}x_{1} & 1\\ \vdots & \vdots \\ x_{i-1} & 1\\ x_{i+1} & 1 \\ \vdots & \vdots\\ x_{n-1} & 1\end{vmatrix}_{(n-2)\,2}\;\;\mbox{for all}\;\; 1 \le i \le n-1
\end{equation}
by Lemma~\ref{lem:AMIRILEMMA32}. Since  $(n-2) + 1$ is even, then Lemma~\ref{lem:NAKAGAMILEMMA20} implies that 
\begin{equation}\label{lem:CullisN31XY:eq106}
\begin{vmatrix}x_{1} & 1\\ \vdots & \vdots \\ x_{i-1} & 1\\ x_{i+1} & 1 \\ \vdots & \vdots\\ x_{n-1} & 1\end{vmatrix}_{(n-2)\,2}
 = 0.
\end{equation}
By combining~\eqref{lem:CullisN31XY:eq105} and~\eqref{lem:CullisN31XY:eq106}, we obtain that
\begin{equation*}
\begin{vmatrix}x_{1} & 1\\ \vdots & \vdots \\ x_{i-1} & 1\\ x_{i+1} & 1 \\ \vdots & \vdots\\ x_{n-1} & 1\\ 0 & 0\end{vmatrix}_{(n-1)\,2} = 0\;\;\mbox{for all}\;\; 1 \le i \le n-1.
\end{equation*}
Thus,
\begin{equation}\label{lem:CullisN31XY:eq103}
\sum_{i=1}^{n-1}(-1)^{i+2}y_i \begin{vmatrix}x_{1} & 1\\ \vdots & \vdots \\ x_{i-1} & 1\\ x_{i+1} & 1 \\ \vdots & \vdots\\ x_{n-1} & 1\\ 0 & 0\end{vmatrix}_{(n-1)\,2} = \sum_{i=1}^{n-1}(-1)^{i+2}y_i \cdot 0 = 0.
\end{equation}

Consider  the second summand on the right-hand side of~\eqref{eq:CullisN31XY:1}. By expanding its second multiplier using Lemma~\ref{lem:NAKAGAMILEMMA20} we obtain that
\begin{equation}\label{lem:CullisN31XY:eq101}
\begin{vmatrix}x_{1} & 1  \\ \vdots & \vdots \\ x_{n-1} & 1\end{vmatrix}_{(n-1)\,2} = \begin{vmatrix}x_{1}\\ \vdots \\ x_{n-1}\end{vmatrix}_{(n-1)\,1}
\end{equation}
because $(n - 1) + 1$ is odd. Using the definition of the Cullis' determinant (Definition~\ref{def:CullisDet}) we conclude that
\begin{equation*}\label{lem:CullisN31XY:eq102}
\begin{vmatrix}x_{1}\\ \vdots \\ x_{n-1}\end{vmatrix}_{(n-1)\,1} =  (-1)^{1}\left(\sum_{1 \le c_1 \le n-1 } (-1)^{c_1} x_{c_1}\right) = \left(\sum_{i=1}^{n-1} (-1)^{i-1} x_{i}\right)
\end{equation*}
and consequently
\begin{equation*}
(-1)^{n+2}y_{n}\begin{vmatrix}x_{1} & 1  \\ \vdots & \vdots \\ x_{n-1} & 1\end{vmatrix}_{(n-1)\,2} = (-1)^{n+1}y_{n}\left(\sum_{i=1}^{n-1} (-1)^{i} x_{i}\right).
\end{equation*}
Since $n$ is odd, then $(-1)^{n+1} = 1$. Hence,
\begin{equation}\label{lem:CullisN31XY:eq104}
(-1)^{n+2}y_{n}\begin{vmatrix}x_{1} & 1  \\ \vdots & \vdots \\ x_{n-1} & 1\end{vmatrix}_{(n-1)\,2} = y_{n}\left(\sum_{i=1}^{n-1} (-1)^{i} x_{i}\right).
\end{equation}

By substituting~\eqref{lem:CullisN31XY:eq105} and~\eqref{lem:CullisN31XY:eq104} instead of the corresponding summands into~\eqref{eq:CullisN31XY:1} we obtain that
\[
\begin{vmatrix}x_{1} & y_{1} & 1\\ \vdots & \vdots & \vdots \\ x_{n-1} & y_{n-1} & 1\\ 0 & y_{n} & 0\end{vmatrix}_{n\,3}
 = 0 + y_n\left(\sum_{i=1}^{n-1} (-1)^{i} x_{i}\right) = y_n\left(\sum_{i=1}^{n-1} (-1)^{i} x_{i}\right)
\]
which establishes the equality~\eqref{lem:CullisN31XY:eq}.
\end{proof}

\begin{lemma}\label{lem:N3LowDegImpliesBInUsupp}Assume that $|\F| > 3$. Let $n \ge 5$ be an odd integer and $B = (b_{i\,j}) \in \M_{n\,3}(\F)$ such that $\rk(B) \ge 2$ and
\begin{equation}\label{lem:N3LowDegImpliesBInU:eq}
   \deg_{\lambda}(\det_{n\,3}(A + \lambda B)) \le 1\;\;\mbox{for all}\;\; A \in \M_{n\,3}(\F).
\end{equation}
\begin{enumerate}[label=(\roman*), ref=(\roman*)]
\item\label{lem:N3LowDegImpliesBInU:part1} If $b_{n\,1} = 0$ and $b_{n\,2} \neq 0$, then $\sum_{i=1}^n (-1)^{i} b_{i\,1} = 0.$
\item\label{lem:N3LowDegImpliesBInU:part2} If there exist $1 \le i_1, i_2 \le n$ such that 
    \begin{equation}\label{lem:N3LowDegImpliesBInU:eq1101}
    b_{i_1\,1} = b_{i_2\,2} = b_{i_1\,3} = 0\;\;\mbox{and}\;\; b_{i_1\,2}, b_{i_2\,1} \neq 0,
    \end{equation} then $B \in U_{n\,3}$.
\end{enumerate}
\end{lemma}
\begin{proof} \textbf{\ref{lem:N3LowDegImpliesBInU:part1}}
  Let $A \in \M_{n\,3}(\F)$ be defined by
  \[
  A = \begin{pmatrix}0 & 0 & 1\\ \vdots & \vdots & \vdots\\0 & 0 & 1\\0 & 0 & 0 \end{pmatrix}.
  \]
  
Let $P \in \F[\lambda]$ be a polynomial defined by $P = \det_{n\, k} (A + \lambda B).$ By Corollary~\ref{cor:DegDetABLEQK}, there exist $a_0, a_1, a_2, a_3 \in \F$ such that
\[
P(\lambda) = a_0 + a_1\lambda + a_2\lambda^2 + a_3 \lambda^3.
\] 

On the one hand, the condition~\eqref{lem:N3LowDegImpliesBInU:eq} means in particular that
\begin{equation}\label{lem:N3LowDegImpliesBInU:eq101}
a_2 = 0.
\end{equation}

On the other hand, \eqref{eq:CullisBinomialExpansion} implies that
\begin{equation}\label{lem:N3LowDegImpliesBInU:eq102}
a_2 = \det_{n\,3} \Bigl(A(|1]\Big| B(|2] \Big| B(|3] \Bigr) +  \det_{n\,3} \Bigl(B(|1]\Big| A(|2] \Big| B(|3] \Bigr)  +  \det_{n\,3} \Bigl(B(|1]\Big| B(|2] \Big| A(|3] \Bigr),
\end{equation}
being a coefficient of $\lambda$ in expansion~\eqref{eq:CullisBinomialExpansion}. Since $A(|1]$ and $A(|2]$ are zero columns, we conclude by splitting off the second and the third term of~\eqref{lem:N3LowDegImpliesBInU:eq102} that

\begin{equation}\label{lem:N3LowDegImpliesBInU:eq103}
a_2 = \det_{n\,3}  \Bigl(B(|1]\Big| B(|2] \Big| A(|3] \Bigr).
\end{equation}

By comparing~\eqref{lem:N3LowDegImpliesBInU:eq101} and~\eqref{lem:N3LowDegImpliesBInU:eq103} we obtain that
\begin{equation}\label{lem:N3LowDegImpliesBInU:eq104}
 \det_{n\,3}  \Bigl(B(|1]\Big| B(|2] \Big| A(|3] \Bigr)  = 0.
\end{equation}

Note that, since $b_{n\, 1} = 0$, then $\det_{n\,3}  \Bigl(B(|1]\Big| B(|2] \Big| A(|3] \Bigr)$ has the same form as the matrix on the left-hand side of~\eqref{lem:CullisN31XY:eq}. Hence, 
\begin{equation}\label{lem:N3LowDegImpliesBInU:eq105}
    \det_{n\,3}  \Bigl(B(|1]\Big| B(|2] \Big| A(|3] \Bigr) = -b_{n\,2}(b_{1\,1} - b_{2\,1} + \ldots - b_{(n-1)\,1})
\end{equation}
      by Lemma~\ref{lem:CullisN31XY}. Substituting the right-hand side of~\eqref{lem:N3LowDegImpliesBInU:eq105} into~\eqref{lem:N3LowDegImpliesBInU:eq104} yields that \[b_{n\,2}\left(b_{1\,1} - b_{2\,1} + \ldots - b_{(n-1)\,1}\right) = 0.\] 
Since $b_{n\,2} \neq 0$ by another assumption of the claim~\ref{lem:N3LowDegImpliesBInU:part1}, \[b_{1\,1} - b_{2\,1} + \ldots - b_{(n-1)\,1} = 0.\] Therefore, 
      \[
        \sum_{i=1}^{n} (-1)^{i}b_{i\,1} = (-1)\bigl((b_{1\,1} - b_{2\,1} + \ldots - b_{(n-1)\,1}) + b_{n\,1}\bigr) = (-1)(0 + 0) = 0.
      \]  

\paragraph{\ref{lem:N3LowDegImpliesBInU:part2}} Let $B^{(1)}, B^{(2)}, B^{(3)} \in \M_{n\,3}(\F)$ be defined by: 
\begin{equation}\label{lem:N3LowDegImpliesBInU:eq11101}
B^{(1)} = \SCS_{(i_1+1)\,1}(B),\,B^{(2)} = \SCS_{(i_2+1)\,2}(B),\,B^{(3)} = \SCS_{(i_1+1)\,3}(B),
\end{equation}
where a linear map $\SCS_{a\,b}$ is defined in Definition~\ref{def:SCSDef}. By Lemma~\ref{lem:ReduceTo11ByShift}, $\SCS_{a\,b}$ preserves $\det_{n\,3}$ by for all $1 \le a \le n, 1 \le b \le 3$. In addition, $\deg_{\lambda}\left(\det_{n\,3}(A + \lambda B)\right) \le 1$ for all $A \in \M_{n\,3}(\F)$ by the conditions of the lemma. Hence, $$\deg_{\lambda}\left(\det_{n\,3}(A + \lambda B^{(j)})\right) \le 1\;\;\mbox{for all}\;\; A \in \M_{n\,3}(\F)\;\; \mbox{and}\;\;1 \le j \le 3$$ by  Lemma~\ref{lem:DetNKPresPresLowDeg}. 

Let us show that
\begin{equation}\label{lem:N3LowDegImpliesBInU:eq11111}
b^{(j)}_{n\,1} = 0\;\;\mbox{and}\;\;b^{(j)}_{n\,2} \neq 0\;\;\mbox{for all}\;\;1 \le j \le 3.
\end{equation}
Indeed, using the definition of $\SCS_{a\,b}$ we obtain that
\begin{gather*}
b^{(1)}_{n\,1} \overset{\eqref{lem:N3LowDegImpliesBInU:eq11101}}{=\joinrel=\joinrel=} -(-1)^{n-(i_1+1)}b_{i_1\,1} \overset{\eqref{lem:N3LowDegImpliesBInU:eq1101}}{=\joinrel=\joinrel=} 0,\\
b^{(1)}_{n\,2} \overset{\eqref{lem:N3LowDegImpliesBInU:eq11101}}{=\joinrel=\joinrel=} -(-1)^{n-(i_1+1)}b_{i_1\,2} \overset{\eqref{lem:N3LowDegImpliesBInU:eq1101}}{=\joinrel\neq\joinrel=} 0,
\end{gather*}
\begin{gather*}
b^{(2)}_{n\,1} \overset{\eqref{lem:N3LowDegImpliesBInU:eq11101}}{=\joinrel=\joinrel=} -(-1)^{n-i_2}b_{i_2\,2} \overset{\eqref{lem:N3LowDegImpliesBInU:eq1101}}{=\joinrel=\joinrel=} 0,\\
b^{(2)}_{n\,2} \overset{\eqref{lem:N3LowDegImpliesBInU:eq11101}}{=\joinrel=\joinrel=} -(-1)^{n-(i_2+1)}b_{i_2\,1} \overset{\eqref{lem:N3LowDegImpliesBInU:eq1101}}{=\joinrel\neq\joinrel=} 0,
\end{gather*}
\begin{gather*}
b^{(3)}_{n\,1} \overset{\eqref{lem:N3LowDegImpliesBInU:eq11101}}{=\joinrel=\joinrel=} -(-1)^{n-i_1}b_{i_1\,3} \overset{\eqref{lem:N3LowDegImpliesBInU:eq1101}}{=\joinrel=\joinrel=} 0\;\;\mbox{and}\\
b^{(3)}_{n\,2}  \overset{\eqref{lem:N3LowDegImpliesBInU:eq11101}}{=\joinrel=\joinrel=} -(-1)^{n-(i_1+1)}b_{i_1\,2} \overset{\eqref{lem:N3LowDegImpliesBInU:eq1101}}{=\joinrel\neq\joinrel=} 0.
\end{gather*}

From~\eqref{lem:N3LowDegImpliesBInU:eq11111} we conclude that each $B^{(j)}$ satisfies the conditions of the claim~\ref{lem:N3LowDegImpliesBInU:part1} and consequently 
\begin{equation}\label{lem:N3LowDegImpliesBInU:eq106}
\sum_{l=1}^n (-1)^{l} b^{(j)}_{l\,1} = 0\;\;\mbox{for all}\;\;1 \le l \le 3.
\end{equation}
By Lemma~\ref{lem:SemiCyclicUKImp} applied to $B$ and each $B^{(j)}$, the implication
\begin{equation}\label{lem:N3LowDegImpliesBInU:eq107}
\sum_{l=1}^n (-1)^{l} b^{(j)}_{l\,1} = 0 \Rightarrow \sum_{l=1}^n (-1)^{l} b_{l\,j} = 0
\end{equation}
holds for every $1 \le j \le 3$. Combining~\eqref{lem:N3LowDegImpliesBInU:eq106} and~\eqref{lem:N3LowDegImpliesBInU:eq107}, we deduce that 
$$\sum_{l=1}^n (-1)^{l} b_{l\,j} = 0\;\;\mbox{for all}\;\;1 \le j \le 3.$$
Thus, $B \in U_{n\,3}$ by the definition of $U_{n\,3}$. 
\end{proof}

\begin{lemma}\label{lem:N3LowDegImpliesBInU}
Assume that $|\F| > 3$. Let $n \ge 5$ be an odd integer and $B = (b_{i\,j}) \in \M_{n\,3}(\F)$. If $\rk(B) \ge 2$ and
\begin{equation*}\label{lem:N3LowDegImpliesBInU:eqq}
   \deg_{\lambda}(\det_{n\,3}(A + \lambda B)) \le 1\;\;\mbox{for all}\;\; A \in \M_{n\,3}(\F),
\end{equation*}
then $B \in U_{n\,3}$.
\end{lemma}
\begin{proof} Let $B = (b_{i\,j}) \in \M_{n\,3}(\F)$ be a matrix satisfying the conditions of the lemma. We perform the proof by transforming the matrix $B$ into the matrix satisfying the conditions of the claim~\ref{lem:N3LowDegImpliesBInU:part2} of Lemma~\ref{lem:N3LowDegImpliesBInUsupp}. For this purpose, we use the appropriate sequence of linear maps preserving $\det_{n\,3}$ and $U_{n\,3}$.

Since $\rk(B) \ge 2$, then there is $1 \le i_1 \neq i_2 \le n$ and $ 1 \le j_1\neq j_2 \le 3$ such that \begin{equation}\label{lem:N3LowDegImpliesBInU:eq10001}
\rk\biggl(B[i_1, i_2|j_1, j_2]\biggr) = 2.
\end{equation}
Given this, we assume without loss of generality that 
\begin{equation}\label{lem:N3LowDegImpliesBInU:eq1002}
b_{i_1\,j_2} \neq 0.
\end{equation}
Since $j_1 \ne j_2$, there exists a permutation $\sigma$ of $\{1, 2, 3\}$ sending $j_1$ to 1 and $j_2$ to $2$. 

Let $M \in \M_{3\,3}(\F)$ be defined by
\[
M = (-1)^{\sgn(\sigma)}\left(E_{1\,\sigma(1)} + E_{2\,\sigma(2)} + E_{3\,\sigma(3)} \right).
\] 
The definition of $M$ implies that 
\begin{equation}\label{lem:N3LowDegImpliesBInU:eq10003}
\det_3(M) = 1.
\end{equation}

Let $B' \in \M_{n\,3}(\F)$ be defined by $B' = BM.$ By direct calculation we obtain that
\[
B'[i_1, i_2|1, 2] = (-1)^{\sgn(\sigma)}B[i_1, i_2|j_1, j_2]
\]
and consequently
\begin{equation}\label{lem:N3LowDegImpliesBInU:eq1005}
\rk\biggl(B'[i_1, i_2|1, 2]\biggr) = \rk\biggl( B[i_1, i_2|j_1, j_2]\biggr) \overset{\eqref{lem:N3LowDegImpliesBInU:eq10001}}{=\joinrel=\joinrel=} 2 .
\end{equation}
In addition, the definition of $B'$ implies that 
\begin{equation}\label{lem:N3LowDegImpliesBInU:eq1001}
b'_{i_1\,2} = (-1)^{\sgn(\sigma)}b_{i_1\,j_2}.
\end{equation}
By combining~\eqref{lem:N3LowDegImpliesBInU:eq1001} and~\eqref{lem:N3LowDegImpliesBInU:eq1002}, we obtain that
\begin{equation}\label{lem:N3LowDegImpliesBInU:eq108}
b'_{i_1\,2}  \neq 0.
\end{equation}

The equality~\eqref{lem:N3LowDegImpliesBInU:eq108} implies that $b'_{i_1\,2}$ is invertible. Hence, the matrix 
$$M' = \begin{pmatrix}1 & 0 & 0\\ -b'_{i_1\,1}\left({b'_{i_1\,2}}\right)^{-1} & 1 & 0 \\ 0 &  0 & 1\end{pmatrix} \in \M_{3\,3}(\F)$$
is properly defined. It follows from the definition of $M'$  that
\begin{equation}\label{lem:N3LowDegImpliesBInU:eq112}
\det_3(M') = 1.
\end{equation}

Let $B'' \in \M_{n\,3}(\F)$ be defined by $B'' = B'M'$. Direct calculations show that
 \begin{equation}\label{lem:N3LowDegImpliesBInU:eq1004}
 b''_{i_1\,1} = 0.
\end{equation}
In addition, since the left multiplication of $B'$ by $M'$ acts as an addition of the second column of $B'$ multiplied by $-b'_{i_1\,1}\left(b'_{i_1\,2}\right)^{-1}$ to the first column of $B'$, then
 \begin{equation}\label{lem:N3LowDegImpliesBInU:eq1006}
 \rk\biggl(B''[i_1, i_2|1, 2]\biggr) = \rk\biggl(B'[i_1, i_2|1, 2]\biggr) \overset{\eqref{lem:N3LowDegImpliesBInU:eq1005}}{=\joinrel=\joinrel=} 2.
 \end{equation}
The equalities~\eqref{lem:N3LowDegImpliesBInU:eq1004} and~\eqref{lem:N3LowDegImpliesBInU:eq1006} together imply that
\begin{equation}\label{lem:N3LowDegImpliesBInU:eq109}
b''_{i_1\,2} \neq 0\;\;\mbox{and}\;\; b''_{i_2\,1} \neq 0,
\end{equation}
because otherwise the matrix $B''[i_1, i_2|1, 2]$ has zero row or zero column and consequently $\rk\biggl(B''[i_1, i_2|1, 2]\biggr) \le 1$.

Thus, $b''_{i_1\,2}$ and $b''_{i_2\,1}$ are invertible by~\eqref{lem:N3LowDegImpliesBInU:eq109}. Hence, the matrix 
$$M'' = \begin{pmatrix}1 &  -b''_{i_2\,2}\left({b''_{i_2\,1}}\right)^{-1} & 0\\ 0 & 1 & -b''_{i_1\,3}\left({b''_{i_1\,2}}\right)^{-1} \\ 0 & 0 & 1\end{pmatrix} \in \M_{3\,3}(\F)$$
is properly defined. It follows from the definition of $M''$ that
\begin{equation}\label{lem:N3LowDegImpliesBInU:eq10004}
\det_3(M'') = 1.
\end{equation}

Let $B''' \in \M_{n\,3}$ be defined by $B''' = B''M''.$
By direct calculations we obtain that 
\begin{equation}\label{lem:N3LowDegImpliesBInU:eq1007}
b'''_{i_1\,3} = b'''_{i_2\,2} = 0,
\end{equation}
\begin{equation}\label{lem:N3LowDegImpliesBInU:eq1008}
b'''_{i_1\,2} = b''_{i_1\,2},\quad b'''_{i_2\,1} = b''_{i_2\,1},
\end{equation}
and
\begin{equation}\label{lem:N3LowDegImpliesBInU:eq1009}
b'''_{i_1\,1} = b''_{i_1\,1}.
\end{equation}
By combining~\eqref{lem:N3LowDegImpliesBInU:eq109} and~\eqref{lem:N3LowDegImpliesBInU:eq1008}, we conclude that
\begin{equation}\label{lem:N3LowDegImpliesBInU:eq1010}
b'''_{i_1\,2} \neq 0\;\;\mbox{and}\;\; b'''_{i_2\,1}\neq 0.
\end{equation}
Similarly,
\begin{equation}\label{lem:N3LowDegImpliesBInU:eq110}
b'''_{i_1\,1} \overset{\eqref{lem:N3LowDegImpliesBInU:eq1009}}{=\joinrel=\joinrel=}  b''_{i_1\,1} \overset{\eqref{lem:N3LowDegImpliesBInU:eq1004}}{=\joinrel=\joinrel=} 0.
\end{equation}

Let $T_M \colon \M_{n\,3}(\F) \to \M_{n\,3}(\F)$ be a linear map defined by  
$$T_M(X)= X\left(MM'M''\right)\;\;\mbox{for all}\;\; X \in \M_{n\,3}(\F).$$ 
Since $\det_3(M) = \det_3(M') = \det_3(M'') = 1$ by~\eqref{lem:N3LowDegImpliesBInU:eq10003}, \eqref{lem:N3LowDegImpliesBInU:eq112}, and~\eqref{lem:N3LowDegImpliesBInU:eq10004} then\linebreak $\det_3(MM'M'') = 1$. Hence, $T_M$ preserves $\det_{n\,3}$ by Lemma~\ref{cor:rightmatrixmult}. The definition of $B'''$ implies that $T_M$ sends $B$ to $B'''$. Hence, 
\begin{equation}\label{lem:N3LowDegImpliesBInU:eq111}
\deg_{\lambda}(\det_{n\,3}(A + \lambda B''')) \le 1\;\;\mbox{for all}\;\;A \in \M_{n\,k}(\F)
\end{equation}
by  Lemma~\ref{lem:DetNKPresPresLowDeg} applied to $T_M$ and $B$.

From~\eqref{lem:N3LowDegImpliesBInU:eq1010}, \eqref{lem:N3LowDegImpliesBInU:eq110} and~\eqref{lem:N3LowDegImpliesBInU:eq111} we obtain that $B'''$ and $i_1, i_2$ satisfy all the conditions of the part \ref{lem:N3LowDegImpliesBInU:part2} of  Lemma~\ref{lem:N3LowDegImpliesBInUsupp} and consequently $B''' \in U_{n\,3}$. 

Since $B''' = B(MM'M'')$ and the matrices $M, M', M''$ are invertible, then\linebreak $B = B'''(MM'M'')^{-1}$. Hence, $B$ is obtained from $B'''$ by elementary column operations and belongs to $U_{n\,3}$ as well as $B'''$. 
\end{proof}
  %

Lemma~\ref{lem:N3LowDegImpliesBInU} allows us to prove that every linear map preserving $\det_{n\,3}$ should preserve the set of matrices having rank one and not belonging to $U_{n\,3}$ (Lemma~\ref{cor:RkAEq1NotInU3ImpliesRkTAEq1}). This fact comprises the significant part of the proof of Lemma~\ref{lem:LinPresRankK3}.

\begin{lemma}\label{cor:RkAEq1NotInU3ImpliesRkTAEq1}Assume that $|\F| > 3$ and  $n \ge 5$ is an odd integer. Let $T\colon \M_{n\,3} (\F) \to \M_{n\,3} (\F)$ be a linear map such that $\det_{n\, 3} (T(X)) = \det_{n\,3}(X)$ for all $X \in \M_{n\,3} (\F)$ and $B \in \M_{n\,3}(\F)$ be a matrix such that $\rk(B) = 1$ and $B \not\in U_{n\,3}$. Then $\rk(T(B)) = 1.$
  \end{lemma}
\begin{proof}
By contradiction. Suppose that $B \in \M_{n\,3}$ is such that $\rk (B) = 1$, 
\begin{equation}\label{cor:RkAEq1NotInU3ImpliesRkTAEq1:eq101}
B \not\in U_{n\,3}
\end{equation} and 
\begin{equation}\label{cor:RkAEq1NotInU3ImpliesRkTAEq1:eq103}
\rk(T(B)) \neq 1.
\end{equation}
Hence, $B \neq O_{n\,3}$ and consequently 
\begin{equation}\label{cor:RkAEq1NotInU3ImpliesRkTAEq1:eq104}
T(B) \neq O_{n\,3}
\end{equation}
because $T$ is invertible by Lemma~\ref{lem:LinPresIso}. It follows from~\eqref{cor:RkAEq1NotInU3ImpliesRkTAEq1:eq103} and~\eqref{cor:RkAEq1NotInU3ImpliesRkTAEq1:eq104} that $\rk(T(B)) \ge 2$. Since $T$ is invertible and $T(U_{n\,3}) = U_{n\,3}$ by Lemma~\ref{lem:DNKPreserverPreservesUn3}, then~\eqref{cor:RkAEq1NotInU3ImpliesRkTAEq1:eq101} implies that
\begin{equation}\label{cor:RkAEq1NotInU3ImpliesRkTAEq1:eq102}
T(B) \not\in U_{n\,3}.
\end{equation}

Our assumption that $\rk (B) = 1$ allows us to conclude by Lemma~\ref{lem:Rank1CullisDeg1} that
$$\deg_{\lambda} (\det_{n\, k} (A + \lambda B)) \le 1\;\;\mbox{for all}\;\;A \in \M_{n\, k}(\F).$$ Therefore,
$$\deg_{\lambda} (\det_{n\, k} (A + \lambda T(B))) \le 1\;\;\mbox{for all}\;\;A \in \M_{n\, k}(\F)$$
by Lemma~\ref{lem:DetNKPresPresLowDeg}. Hence, $T(B) \in U_{n\,3}$ by Lemma~\ref{lem:N3LowDegImpliesBInU}. This contradicts with~\eqref{cor:RkAEq1NotInU3ImpliesRkTAEq1:eq102}.
\end{proof}

\begin{lemma}\label{lem:LinPresRankK3}Assume that $|\F| > 3$ and  $n \ge 5$ is an odd integer. Let $T\colon \M_{n\, 3} (\F) \to \M_{n\, 3} (\F)$ be a linear map such that $\det_{n\, 3} (T(X)) = \det_{n\,3}(X)$ for all $X \in \M_{n\, 3} (\F).$ Then $\rk (X) = 1$ implies $\rk (T(X)) =  1$ whenever $X \in \M_{n\, 3}(\F)$.
\end{lemma}

\begin{proof}Let $X \in \M_{n\,3}(\F)$ be such that $\rk (X) = 1$. Hence, there exists $\mathbf 0 \neq \mathbf v \in \F^3$ and $\lambda_1, \ldots, \lambda_n \in \F$ such that $(\lambda_1, \ldots, \lambda_n) \neq \mathbf 0$ and
\begin{equation}\label{lem:LinPresRankK3:eq1001}
(x_{i\,1}, x_{i\,2}, x_{i\,3}) = \lambda_i\mathbf v^t\;\;\mbox{for all}\;\;1 \le i \le n.
\end{equation}
Let $1 \le j_0 \le 3$ be the index of the nonzero coordinate of $\mathbf v$, i.e. 
\begin{equation}\label{lem:LinPresRankK3:eq100002}
\mathbf v_{j_0} \neq 0.
\end{equation}

Let $B^{(1)}, \ldots, B^{(n)} \in \M_{n\,3}(\F)$ be defined by
\[
B^{(i)} = \mathbf v_1 E_{i\,1} + \mathbf v_2 E_{i\,2} + \mathbf v_3 E_{i\,3}\;\;\mbox{for all}\;\;1 \le i \le n.
\]
That is,
\[
B^{(1)} = \begin{pmatrix}\mathbf v_1 & \mathbf v_2 &  \mathbf v_3\\ \vdots & \vdots & \vdots\\ 0 & 0 &  0 \end{pmatrix},\quad  \ldots, \quad B^{(n)} = \begin{pmatrix} 0 & 0 & 0\\ \vdots & \vdots & \vdots\\ \mathbf v_1 & \mathbf v_2 &  \mathbf v_3 \end{pmatrix}.
\]
From~\eqref{lem:LinPresRankK3:eq1001} we conclude that
\begin{equation}\label{lem:LinPresRankK3:eq1002}
T(X) = \lambda_1 T\left(B^{(1)}\right) + \ldots +  \lambda_n T\left(B^{(n)}\right).
\end{equation}

Let us show that the sequence $c_1 = -1, c_2 = 1, \ldots, c_n = (-1)^n \in \F$ and the matrices $T\left(B^{(1)}\right),$ $\ldots,$ $T\left(B^{(n)}\right)$  satisfy the conditions of Lemma~\ref{lem:Rk1NonzeroRk1All}. 


Let $\alpha_1, \ldots, \alpha_n \in \F$ be such that 
\begin{equation}\label{lem:LinPresRankK3:eq101}
c_1 \alpha_1 + \ldots + c_n \alpha_n = \sum_{i=1}^{n}(-1)^{i} \alpha_i \neq 0.
\end{equation}
Let $B = (b_{i\,j})  \in \M_{n\,3}(\F)$ be defined by 
\[
B = \alpha_1 B^{(1)} + \ldots + \alpha_n B^{(n)}.
\]

Since 
\[
\sum_{i=1}^n (-1)^{i}b_{i\, j_0} = \sum_{i=1}^n (-1)^{i}\alpha_i \mathbf v_{j_0} = \left(\sum_{i=1}^n (-1)^{i}\alpha_i \right)\mathbf v_{j_0}
\] 
by the definition of $B$ and $B^{(1)}, \ldots, B^{(n)}$, then 
\begin{equation}\label{lem:LinPresRankK3:eq103}
\sum_{i=1}^n (-1)^{i}b_{i\, j_0} = \left(\sum_{i=1}^n (-1)^{i}\alpha_i \right)\mathbf v_{j_0} \neq 0
\end{equation}
by~\eqref{lem:LinPresRankK3:eq100002} and~\eqref{lem:LinPresRankK3:eq101}. Hence, 
\begin{equation}\label{lem:LinPresRankK3:eq104}
B = \alpha_1 B^{(1)} + \ldots + \alpha_n B^{(n)}\not\in U_{n\, 3}
\end{equation}
by~\eqref{lem:LinPresRankK3:eq103}. 
Since $T$ is invertible by Lemma~\ref{lem:LinPresIso} and $T(U_{n\,3}) = U_{n\,3}$ by Lemma~\ref{lem:DNKPreserverPreservesUn3}, then~\eqref{lem:LinPresRankK3:eq104} implies that
\begin{equation*}
 T(B) = \alpha_1 T(B^{(1)}) + \ldots + \alpha_3 T(B^{(n)})\not\in U_{n\, 3}.
\end{equation*}
The definition of $B$ implies that $\rk(B) = 1$. Therefore, 
$$\rk(T(B)) = \rk\left(\alpha_1 B^{(1)} + \ldots + \alpha_3 B^{(n)}\right) = 1$$
by Lemma~\ref{cor:RkAEq1NotInU3ImpliesRkTAEq1}. 

Thus, the sequence $c_1 = -1, c_2 = 1, \ldots, c_n = (-1)^n \in \F$ and the matrices $T\left(B^{(1)}\right),$ $\ldots,$ $T\left(B^{(n)}\right)$  satisfy the conditions of Lemma~\ref{lem:Rk1NonzeroRk1All}. Since 
$$T(X) = \lambda_1 T\left(B^{(1)}\right) + \ldots +  \lambda_n T\left(B^{(n)}\right)$$
 by~\eqref{lem:LinPresRankK3:eq1002},  we conclude from Lemma~\ref{lem:Rk1NonzeroRk1All} that 
\begin{equation}\label{lem:LinPresRankK3:eq102}
\rk(T(X)) \le 1.
\end{equation}
Recall that $\rk(X) = 1$ by the initial assumption. Hence, $X \neq O_{n\,3}$. The linear map $T$ is invertible by Lemma~\ref{lem:LinPresIso}. Therefore,
\begin{equation}\label{lem:LinPresRankK3:eq100001}
T(X) \neq O_{n\,3}.
\end{equation}
By combining~\eqref{lem:LinPresRankK3:eq102} and~\eqref{lem:LinPresRankK3:eq100001} we conclude that $\rk(T(X)) = 1$.
\end{proof}

\begin{theorem}\label{thm:MainTheoremEvenINCLK3}Assume that  $|\F| > 3$, $n \ge 5$ and $n$ is odd. Suppose that $T\colon \M_{n\, 3} (\F) \to \M_{n\, 3} (\F)$ is a linear map. Then $\det_{n\, 3} (T(X)) = \det_{n\,3}(X)$ for all $X \in \M_{n\, 3} (\F)$ if and only if there exist $A \in \M_{n\, n}(\F)$ and $B \in \M_{3\, 3}(\F)$ such that
\begin{equation}\label{thm:MainTheoremEvenINCLK3:eq}
\det_{n\, 3} \Bigl(A[i_1,i_2, i_3|)\Bigr) \det_3 (B) = (-1)^{i_1 + i_2 + i_3 - 1 - 2 - 3}
\end{equation}
for all increasing sequences $1 \le i_1 < i_2 < i_3 \le n $ and
\[
T(X) = AXB\;\;\mbox{for all}\;\;X \in \M_{n\, 3} (\F).
\]
\end{theorem}

\begin{proof}The sufficiency follows from Lemma~\ref{lem:TwoSidedMulPreservesDet}. Let us prove the necessity.

It follows from Lemma~\ref{lem:LinPresRankK3} that $\rk (T(X)) = 1$ for all $X$ with $\rk (X) = 1.$ In addition, $T$ is invertible by Lemma~\ref{lem:LinPresIso}. Therefore, by Lemma~\ref{lem:westwick} there exist invertible matrices $A \in \M_{n\,n}(\F)$ and $B \in \M_{3\,3}(\F)$ such that
  \[
    T(X) = AXB\;\;\mbox{for all}\;\;X \in \M_{n\,3}(\F).
  \]
 These $A$ and $B$ satisfy the condition~\eqref{thm:MainTheoremEvenINCLK3:eq} by Lemma~\ref{lem:TwoSidedMulPreservesDet} because $T$ preserves $\det_{n\,3}$.
\end{proof}

In order to conclude Theorem~\ref{thm:MainTheoremOddINCLK3} from Theorem~\ref{thm:MainTheoremEvenINCLK3}, we  use the same approach as in~\cite{Guterman2025b}. This approach relied on the definitions and propositions provided below.

\begin{definition}[{\cite[Definition 3.6]{Guterman2025b}}]\label{def:WNK}By $W_{n\,k}\subseteq \M_{n\,k}(\mathbb F)$ we denote the vector space of matrices whose rows are all equal. That is,
\[
W_{n\,k} = \{\begin{pmatrix}y_1 & \cdots & y_k\\ \vdots & \ddots & \vdots\\ y_1 & \cdots & y_k\end{pmatrix} \mid y_1, \ldots, y_k \in \F\}.
 \]
\end{definition}

\begin{lemma}[{\cite[Corollary~2.18]{Guterman2025b}}]\label{lem:AdditionLinearMapNPLusKOdd} Assume that $n + k$ is odd. Let 
$\phi \colon \M_{n\,k}(\F) \to W_{n\,k}$ and $T \colon \M_{n\,k}(\F) \to \M_{n\,k}(\F)$ be linear maps such that  $\det_{n\, k} (T(X)) = \det_{n\,k}(X)$ for all $X \in \M_{n\,k}(\F)$. Denote by $S$ a linear map $\M_{n\,k}(\F) \to \M_{n\,k}(\F)$ defined by
\begin{equation}\label{eq:AddingKernel}
S(X) = T(X) + \phi (X)\;\;\mbox{for all}\;\;X \in \M_{n\,k}(\F).
\end{equation}
 Then
\[
\det_{n\, k} (S(X)) = \det_{n\,k}(X)\;\;\mbox{for all}\;\;X \in \M_{n\,k} (\F).
\]
\end{lemma}

\begin{definition}[{\cite[Definition~4.1]{Guterman2025b}}]
  By $\M^{0}_{n\,k}(\F) \subseteq \M_{n\,k}(\F)$ we denote the vector space of all matrices such that all entries in their last row are equal to zero, i.e.
  \[
    \M^{0}_{n\,k}(\F) = \{X \in \M_{n\,k}(\F) \mid x_{n\,j} = 0, 1 \le j \le k\}.
  \]
\end{definition}

\begin{definition}[{\cite[Definition~4.2]{Guterman2025b}}]\label{def:Lminus}
For every matrix $X \in \M_{n\,k}(\mathbb F)$, by $L^{-}(X) \in \M_{n-1\,k} (\mathbb F)$ we denote a matrix defined by
\[
L^{-}(X) = \begin{pmatrix}X[1|) - X[n|)\\ \vdots \\X[n-1|) - X[n|)\end{pmatrix}.
\]

\end{definition}

\begin{claim}[{\cite[Lemma~4.3]{Guterman2025b}}]\label{claim:LminusMatrixForm}$L^{-}(X) = M_{L^{-}}X$ for all $X \in \M_{n\,k}(\F)$, where
\[
M_{L^{-}} =
\begin{pmatrix}
1 & 0 & \ldots & 0 & -1\\
0 & 1 & \ldots & 0 & -1\\
\vdots & \vdots & \ddots & \vdots & \vdots\\
0 & 0 & \ldots & 1 & -1
\end{pmatrix}
\in \M_{n-1\,n}(\F).
\]
  \end{claim}

\begin{definition}[{\cite[Definition~4.5]{Guterman2025b}}]\label{def:Lplus}
For every matrix $X \in \M_{(n-1)\,k}(\mathbb F)$, by $L^{+}(X) \in \M_{n\,k}(\F)$ we denote the matrix obtained from $X$ by adjoining to it the zero row as the last row; that is
\[
L^{+}(X) = \begin{pmatrix}X\\ \mathbf{0}
\end{pmatrix}.
\]
\end{definition}

\begin{claim}[{\cite[Lemma~4.6]{Guterman2025b}}]\label{claim:LplusMatrixForm}$L^{+}(X) = M_{L^{+}}X$ for all $X \in \M_{(n-1)\,k}(\F)$, where
\[
M_{L^{+}} =
\begin{pmatrix}
1 & 0 & \ldots & 0\\
0 & 1 & \ldots & 0\\
\vdots & \vdots & \ddots & \vdots\\
0 & 0 & \ldots & 1\\
0 & 0 & \ldots & 0
\end{pmatrix}
\in \M_{n\,(n-1)}(\F).
\]
  \end{claim}

\begin{lemma}[{Cf.~\cite[Lemma~4.8]{Guterman2025b}}]\label{lem:LiftDownDetPres}Assume that  $|\F| > k$, $n > k$ and $n + k$ is odd. Let $T$ be a linear map on $\M_{n\, k} (\F)$ such that $\det_{n\, k} (T(X)) = \det_{n\,k}(X)$ for all $X \in \M_{n\,k}(\F)$. If a linear map $S$ on $\M_{(n-1)\,k}(\F)$ is defined by 
$$S = L^{-} \circ T \circ L^{+},$$
then $\det_{(n-1)\, k} (S(Y)) = \det_{(n-1)\,k}(Y)$ for all $Y \in \M_{(n-1)\,k}(\F)$.
\end{lemma}

\begin{lemma}[{Cf.~\cite[Lemma~4.10]{Guterman2025b}}]\label{lem:FromNKOddToNKEvenIfKerW}
Assume that  $|\F| > k$, $n > k$ and $n + k$ is odd. Let $T\colon \M_{n\, k} (\F) \to \M_{n\, k} (\F)$ be a linear map such that $\det_{n\, k} (T(X)) = \det_{n\,k}(X)$ and $\Img(T) \subseteq \M^0_{n\,k}(\F)$. If a linear map $S\colon \M_{(n-1)\,k}(\F) \to \M_{(n-1)\,k}(\F)$ is defined by 
\begin{equation}\label{lem:FromNKOddToNKEvenIfKerW:eq}
S = L^{-} \circ T \circ L^{+},
\end{equation}
then $T = L^{+} \circ S \circ L^{-}$.
\end{lemma}

\begin{lemma}[{Cf.~\cite[Lemma~4.11]{Guterman2025b}}]\label{lem:MatRepLift}Assume that  $|\F| > k$, $n > k$ and $n + k$ is odd. Let $S$ be a linear map on $\M_{(n-1)\,k}(\F)$ and $T$ be a linear map on $\M_{n\,k}(\F)$ defined by 
\[
T = L^{+}\circ S \circ L^{-}.
\] If there exist $A' \in \M_{(n-1)\,(n-1)}(\F)$ and $B' \in \M_{k\,k}(\F)$ satisfying the condition~\eqref{eq:TwoSideMul} such that 
\[
S(Y) = A'YB'\;\;\mbox{for all}\;\; Y \in \M_{(n-1)\,k}(\F),
\]
then there exist $A \in \M_{n\, n}(\F)$ and $B \in \M_{k\, k}(\F)$ satisfying the condition~\eqref{eq:TwoSideMul} such that
\[
T(X) = AXB\;\;\mbox{for all}\;\; X \in \M_{n\,k}(\F).
\]
\end{lemma}

\begin{theorem}\label{thm:MainTheoremOddINCLK3}Assume that $|\F| > 3, n \ge 5$ and $n$ is even. Let $T\colon \M_{n\,3} (\F) \to \M_{n\,3} (\F)$ be a linear map. Then $\det_{n\, 3} (T(X)) = \det_{n\,3}(X)$ for all $X \in \M_{n\,3} (\F)$ if and only if  there exist $A \in \M_{n\, n}(\F)$ and $B \in \M_{3\, 3}(\F)$ such that
\begin{equation}\label{thm:MainTheoremCullisNKOdd:eqq}
\det_{n\, 3} \Bigl(A[i_1,i_2, i_3|)\Bigr) \det_3 (B) = (-1)^{i_1 + i_2 + i_3 - 1 - 2 - 3}
\end{equation}
for all increasing sequences $1 \le i_1 < i_2 < i_3 \le n $ and a linear map $\phi\colon \M_{n\,3}(\F) \to W_{n\,3}$ such that
\begin{equation}\label{thm:MainTheoremCullisNKOdd:eq}
T(X) = AXB + \phi(X)\;\;\mbox{for all}\;\; X \in \M_{n\,3} (\F).
\end{equation}
\end{theorem}

\begin{proof}
\textbf{1.} Let us prove the sufficiency. Let $T$ be a linear map which satisfies the condition~\eqref{thm:MainTheoremCullisNKOdd:eq} for some $A \in \M_{n\, n}(\F)$, $B \in \M_{3\, 3}(\F)$ satisfying the condition~\eqref{thm:MainTheoremCullisNKOdd:eqq} and a linear map $\phi\colon \M_{n\,3}(\F) \to W_{n\,3}$. 

Lemma~\ref{lem:TwoSidedMulPreservesDet} implies that $\left(T - \phi\right)$ preserves $\det_{n\,3}$. Therefore, $T = \left(T - \phi\right) + \phi$ preserves $\det_{n\,3}$ by Lemma~\ref{lem:AdditionLinearMapNPLusKOdd} applied to $(T - \phi)$ and $\phi$. 

\textbf{2.} Let us prove the necessity. Assume that $T\colon \M_{n\,3} (\F) \to \M_{n\,3} (\F)$ is a linear map preserving $\det_{n\, 3}$. Let $\phi\colon \M_{n\,3}(\F) \to W_{n\,3}$ be a linear map defined by
\[
\phi(X) = \begin{pmatrix}T(X)[n|)\\\vdots\\ T(X)[n|)\end{pmatrix},
\]
$\tilde{T}\colon \M_{n\,3} (\F) \to \M_{n\,3} (\F)$ be a linear map defined by $\tilde{T} = T - \phi$ and $S$ be a linear map on $\M_{(n-1)\,3}(\F)$  defined by $S = L^{-}\circ \tilde{T} \circ L^{+}$. 

The definition of $\tilde{T}$ implies that $\Img(\tilde{T}) \subseteq \M^{0}_{n\,3}(\F)$. In addition, $\tilde{T}$ preserves $\det_{n\,3}$ by Lemma~\ref{lem:AdditionLinearMapNPLusKOdd} applied to $T$ and $(-\phi)$. Thus, $\tilde{T}$ satisfies the conditions of Lemma~\ref{lem:FromNKOddToNKEvenIfKerW} and consequently 
\begin{equation}\label{thm:MainTheoremCullisNKOdd:eq1}
\tilde{T} = L^{+}\circ S \circ L^{-}.
\end{equation}

Since $\tilde{T}$ preserves $\det_{n\,3}$, then $S$ preserves $\det_{(n-1)\,3}$ by Lemma~\ref{lem:LiftDownDetPres}.\linebreak Theorem~\ref{thm:MainTheoremEvenINCLK3} implies that there exist $A' \in \M_{(n-1)\,(n-1)}(\F)$ and $B' \in \M_{3\,3}(\F)$ satisfying the condition~\eqref{thm:MainTheoremEvenINCLK3:eq} such that 
\begin{equation}\label{thm:MainTheoremCullisNKOdd:eq2}
S(Y) = A'YB'\;\;\mbox{for all}\;\; Y \in \M_{(n-1)\,3}(\F).
\end{equation}

Thus, the equalities~\eqref{thm:MainTheoremCullisNKOdd:eq1} and~\eqref{thm:MainTheoremCullisNKOdd:eq2} allow us to conclude from Lemma~\ref{lem:MatRepLift} that there exist $A \in \M_{n\, n}(\F)$ and $B \in \M_{3\, 3}(\F)$ satisfying the condition~\eqref{thm:MainTheoremEvenINCLK3:eq} such that 
$$\tilde{T}(X) = AXB\;\;\mbox{for all}\;\; X \in \M_{n\,3}(\F).$$
Therefore,
\[
T(X) = \tilde{T}(X) + \phi(X) = AXB + \phi(X)\;\;\mbox{for all}\;\;X \in \M_{n\,3} (\F).
\]
\end{proof}
\begin{remark}One may see that the only difference of this proof from the proof of Theorem~\ref{thm:THMNPLUSKODD} is that Theorem~\ref{thm:MainTheoremEvenINCLK3} is used instead of Theorem~\ref{thm:THMNPLUSKEVEN}.
\end{remark}

\begin{remark}Theorem~\ref{thm:MainTheoremEvenINCLK3} and Theorem~\ref{thm:MainTheoremOddINCLK3} imply that the linear maps preserving $\det_{n\,3}$ for $n \ge 5$ have the same description as linear maps preserving $\det_{n\,k}$ for $k \ge 4$ and $n \ge k + 2$.
\end{remark}

\bibliographystyle{specialmatrices}

\end{document}